\documentclass{jfm}
\usepackage{color}
\usepackage{amstext}
\usepackage{amssymb}
\usepackage{graphicx}
\usepackage{amsmath}
\usepackage{tcolorbox}
\usepackage{nccmath}
\usepackage{amsfonts}
\usepackage[utf8]{inputenc}

\DeclareMathAlphabet\mathbfcal{OMS}{cmsy}{b}{n}

\newcommand\const{\mathrm{const}}

\newcommand\vf{\boldsymbol{f}}

\newcommand\vn{\boldsymbol{n}}
\newcommand\vu{\boldsymbol{u}}
\newcommand\vv{\boldsymbol{v}}

\newcommand\vx{\boldsymbol{x}}
\newcommand\vy{\boldsymbol{y}}

\begin{document}

{\title[Distinguished Limits and Drifts: between Nonuniqueness and Universality] {Distinguished Limits and Drifts:\\ between Nonuniqueness and Universality}}

\author[V. A. Vladimirov]
{V.\ns A.\ns V\ls l\ls a\ls d\ls i\ls m\ls i\ls r\ls o\ls v}

\affiliation{York University and Leeds University, UK; vv500@york.ac.uk}

%\pubyear{2014} \volume{xx} \pagerange{xx-xx}
\date{August 1st, 2021, York, Ferriby}

\setcounter{page}{1}\maketitle \thispagestyle{empty}

\begin{abstract}

This paper deals with a version of the two-timing method which describes various `slow' effects caused by externally imposed `fast' oscillations.
Such small oscillations are often called \emph{vibrations} and the research area can be referred as \emph{vibrodynamics}.
The governing equations represent a generic system of first-order ODEs containing a prescribed oscillating velocity $\vu$, given in a general form.
Two basic small parameters stand in for the inverse frequency and the ratio of two time-scales; they appear in equations as regular perturbations.
The proper connections between these parameters yield the \emph{distinguished limits}, leading to the existence of closed systems of asymptotic equations.
The aim of this paper is twofold: (i) to clarify (or to demystify) the choices of a slow variable, and
(ii) to give a coherent exposition which is accessible for practical users in applied mathematics, sciences and engineering.
We focus our study on the usually hidden aspects of the two-timing method such as the  \emph{uniqueness or multiplicity of distinguished limits} and  \emph{universal structures of averaged equations}.
The main result is the demonstration that there are two (and only two) different distinguished limits.
The explicit instruction for practically solving ODEs for different classes of $\vu$ is presented.
The key roles of drift velocity and the qualitatively new appearance of the linearized equations are discussed.
To illustrate the broadness of our approach, two examples from mathematical biology are shown.

\medskip

\noindent\textsc{R\'esum\'e.}

Cet article traite d'une version de la méthode à deux temps qui décrit divers effets "lents" causés par des oscillations "rapides" imposées de l'extérieur. 
Ces petites oscillations sont souvent appelées vibrations et le domaine de recherche peut être appelé \emph{vibrodynamique}. 
Les équations gouvernantes considérées représentent un système générique d'EDOs du premier ordre, contenant une vitesse d'oscillation prescrite $\vu$, donnée sous une forme générale. 
Deux petits paramètres de base représentent la fréquence inverse et le rapport de deux échelles de temps ; ils apparaissent dans les équations sous forme de perturbations régulières. 
Les connexions appropriées entre ces paramètres donnent les \emph{limites distinguées}, menant à l'existence de systèmes fermés d'équations asymptotiques. 
L'objectif de cet article est double : (i) clarifier les choix de la variable lente, et (ii) donner une exposition cohérente de la méthode, appliquée à la vibrodynamique. 
Nous concentrons notre étude sur les aspects habituellement cachés de la méthode à deux temps, tels que \emph{l'unicité ou la multiplicité des limites distinguées et les structures universelles des équations moyennées}. Le résultat principal est la démonstration qu'il existe deux (et seulement deux) limites distinguées différentes. L'instruction explicite pour la manipulation pratique des EDOs pour différentes classes de $\vu$ est présentée. Les rôles clés de la vitesse de dérive sont discutés. Pour illustrer l'étendue de notre approche, deux exemples issus de la biologie mathématique sont présentés. Cet article est accessible aux étudiants en mathématiques appliquées, en physique et en génie.


\end{abstract}

\section{Introduction \label{sect01}}

The two-timing method is a classical tool of applied mathematics and perturbation methods, see \emph{e.g.} \cite{Bog, Nayfeh, Grimshaw, Hinch, Kevorkian, Verh, Arn},\newline \cite{San}.
The literature on its engineering applications is also very extensive, see  \emph{e.g.} \cite{Blekhman1, Thom, Blekhman2}.
There are two main classes of problems treated by the two-timing method.
The first one deals with the evolution of  modulated oscillation, resulting in slow changes of various periodic motions;
a typical example is a nonlinear interaction between two pendulums, see \emph{e.g.}  \cite{Krylov, OMalley}.
The second class exposes the `slowly developing' large amplitude effects caused by the imposed `fast' oscillations of small amplitudes.
Such oscillations are known as \emph{vibrations}, the related research area deserves a special name of \emph{vibrodynamics}, see \emph{e.g.} \cite{Vladimirov, Yudovich, Vladimirov1}.
Its paradigm is the Stephenson-Kapitza pendulum also known as the inverted pendulum, see \cite{Kapitza1, Kapitza2}.
In this paper we consider a system of first-order ODEs in variable $t$ that may be seen as a motion $\vx(t)$ of a material particle in a given oscillating velocity field $\vu(\vx,t)$.
The essence of the method is replacing the single variable $t$ with two independent variables: \emph{slow time} $s$ and \emph{fast time} $\tau$.
It leads to an \emph{auxiliary PDE}, that should be solved with the use of standard asymptotic methods.
Then, returning back to the original variable $t$ produces solutions of the original ODE.

The objectives of our  research are two general (but usually hidden) issues, related to the auxiliary PDEs:
(i) \emph{uniqueness and/or multiplicity of distinguished limits} (abbreviated as \emph{DL}s), and (ii) \emph{universality in the structure of averaged equations}.
There are different definitions of \emph{distinguished limits} for PDEs, see, \emph{e.g.} \cite{Nayfeh, Kevorkian}, in which they are used for singular perturbations only.
Some other papers mention \emph{DL}s  in a more general context, see \emph{e.g.} \cite{Klein}.
All definitions operate with `proper relations' between different terms in the dimensionless governing equations, and such relations must lead to the closed systems of asymptotic equations.
For example, if an equation possesses two independent small parameters $\mu$ and $\varepsilon$, then one should define what is the interrelation between their magnitudes.
Otherwise, any successive approximations cannot work.
To resolve this issue, one must study all possible paths in the plane $(\mu,\varepsilon)$, such that $(\mu,\varepsilon)\to (0,0)$.
They may be represented as  $\mu(\varepsilon)$.
If $\varepsilon\to 0$ leads to a \emph{closed system of asymptotic equations} then this path is called \emph{a distinguished limit (DL)}.
A commonly confusing feature is: \emph{DLs} are usually used as \emph{`know-how'}, without any comments to the method of their finding and to the presence or absence of any alternatives.
Most often, \emph{DLs} are used implicitly, without even mentioning them.
The existence of \emph{DL}s demonstrates that the structure of an equation dictates the proper placing of a small parameter within it.
The terms \emph{asymptotic equation/solution/approximation} can also have different meanings.
Usually, only the existence of the first/main term in the related series in $\varepsilon$ is considered, while no attention is paid to the next approximations, see \emph{e.g.} \cite{Kevorkian}.
An opposite limiting case is a \emph{regular asymptotic procedure} that implies that the successive approximation of any order in  $\varepsilon$ can be calculated or at least the related closed systems of equations can be derived.

An important issue is the mathematical justification of the approximations obtained by the two-timing method.
Such studies employ functional analysis and deal with the upper bounds of the differences between exact and asymptotic solutions.
This research area is still actively developing with a number of different results and authors, see \emph{e.g.} \cite{Simonenko, Murdock0, Murdock, Yudovich}\newline\cite{San, Leven, Leven-a}.
%,\newline \cite{San}.
Some of them prove that asymptotic solutions are valid in the \emph{expanding time-interval} of length, for instance $O(1/\varepsilon)$, whilst some others derive results for any finite interval $(0,\mathcal{T}\,)$,
where multipliers in the obtained upper bounds may increase with $\mathcal{T}$, even exponentially.
Such results for a class of first-order ODEs and any finite time-interval $(0,\mathcal{T}\,)$ can be found in \cite{Leven, Leven-a}.
At the same time, applied users usually believe that there is significant experimental and computational evidence that asymptotic solutions often work far beyond the interval of their mathematical approval, see \emph{e.g.} \cite{Thom, Yudovich, Blekhman2}.
In this paper we adopt \emph{a physical style of exposition}, avoiding the pure mathematical issues.
In particular, we look for solutions in the form of a regular series or polynomial in $\varepsilon$, but every time we take notice how many terms we have explicitly calculated.

We start the paper by obtaining the dimensionless form of the equations.
This step is often underestimated; different scalings can produce the different appearances of  small parameters in the equation and different asymptotic models.
Then we analyse the \emph{multiplicity and uniqueness  of} \emph{DL}s and  expose the \emph{universal structures of averaged equations}.
In doing this, we describe and use the \emph{DL-procedure}.
As the main result, we obtain two and only two \emph{DL}s (called \emph{DL-1} and \emph{DL-2}), where both may lead to \emph{regular asymptotic solutions} represented by infinite series.
Our explicit derivations are restricted by first three successive approximations that lead to the closed systems of equations.
It is remarkable that these equations possess \emph{universal structures}; they consist of (i) the amplitude decomposition of equations obtained in the absence of vibrations and (ii) a quadratic (in the oscillating part of velocity $\widetilde{\vu}$) averaged function called the \emph{drift velocity} ${\overline{\mathbfcal{V}}_2}$, which we require to be nonzero.
\emph{DL-2} leads to a closed averaged system of equations in the zeroth approximation, which  include ${\overline{\mathbfcal{V}}_2}$.
Obtaining equations of such type always represents the main aim of applied researchers.
In \emph{DL-1}, the averaged equations of the zeroth approximation formally coincide with the original system.
The first \emph{DL-1} approximation exhibits a similarity with a linearized version of zero-order equations, however it is complimented by ${\overline{\mathbfcal{V}}_2}$ as a nonhomogeneous `driving' term;
this reveals a qualitatively new impact of linear approximation.
The drift velocity ${\overline{\mathbfcal{V}}_2}$ generalises the classical notion of  \emph{drift} in plasma physics and in fluid dynamics, see \emph{e.g.} \cite{Stokes, Lamb, Batchelor, Craik, Eames, Moni}.
In addition, we use a relatively new general term in our explanations: various averaged products of oscillating functions are called \emph{vibrogenic terms}.
Then the drift velocity ${\overline{\mathbfcal{V}}_2}$ may  also be  called \emph{the vibrogenic velocity}.
Such an addition in terminology was introduced by \cite{Yudovich} to make the explicit distinction between oscillating functions and their averaged nonlinear counterparts.

The structure of the paper is the following:
two sections, \emph{Setting of two-timing problem} and \emph{Distinguished limit procedure}, describe the approach used.
The next section \emph{Derivations} contains all the required calculations.
Two sections, \emph{Summary of results} and \emph{Practical instruction}, contain the formulation of results obtained and the explicit instruction for practical handling of ODEs in vibrodynamics, for different classes of $\vu$.
Interpretations and extensions of our results are given in the sections \emph{Drift velocity} and \emph{Discussion}.
For the sake of brevity, our exposition uses terminology from classical mechanics and fluid dynamics.
However, the described approach may be useful for studying applied differential equations from various research areas, after introducing time-oscillations/vibrations to their coefficients.
Two examples in the most renowned models of mathematical biology are given in the subsection \emph{Examples of drifts}.

\section{Setting of two-timing problem}%\label{sect1}

Let us consider a system of  first-order ODEs describing a motion $\vx^\dag(t^\dag)$ of a particle with a prescribed velocity field (or in a prescribed flow) $\vu^\dag(\vx^\dag, t^\dag)$
\begin{eqnarray}\label{TT-ode}
 d\vx^\dag/dt^\dag=\vu^\dag(\vx^\dag, t^\dag),
\end{eqnarray}
where $\vx^\dag=(x_1^\dag,x_2^\dag,x_3^\dag)$ is a cartesian coordinate,  $t^\dag$ is time,
$\vu^\dag(\vx^\dag, t^\dag)$ is a given function, and daggers denote dimensional variables.
The spatial domain $\vx\in\mathcal{D}^\dag$ can be finite or infinite, it does not effect our consideration.
We consider oscillatory motions that possess characteristic scales of velocity and length, together with two additional time-scales:
\begin{eqnarray}\label{scales}
U^\dag, \quad  L^\dag, \quad  T^\dag_\text{slow},\quad T^\dag_\text{fast}.
\end{eqnarray}
There are therefore two independent dimensionless parameters
\begin{eqnarray}\label{parameters}
T_\text{slow}\equiv  T^\dag_\text{slow}/T^\dag,  \quad T_\text{fast}\equiv T^\dag_\text{fast}/T^\dag,\quad\text{where}\quad T^\dag\equiv L^\dag/U^\dag,
\end{eqnarray}
which represent the dimensionless time-scales.
The scale $T_\text{fast}^\dag$ characterises the given period of oscillations; hence the dimensionless and dimensional frequencies of oscillations are:
\begin{eqnarray}\label{frequency}
\omega^\dag \equiv 1/T^\dag_\text{fast}, \quad \omega\equiv T^\dag/T^\dag_\text{fast}.
\end{eqnarray}
The scale $T_\text{slow}^\dag$ remains undefined, it will be calculated in the process of finding the distinguished limits (\emph{DL}s).
We choose the dimensionless  variables as
\begin{eqnarray}\label{x-t}
\vx\equiv\vx^\dag/L^\dag, \quad t\equiv t^\dag/T^\dag.
\end{eqnarray}
The dimensionless `fast time' $\tau$ and `slow time' $s$ are defined as:
\begin{eqnarray}\label{tay-s}
\tau \equiv t/T_\text{fast}\equiv\omega t, \quad s\equiv t/T_\text{slow}\equiv St,\quad \text{with}\quad S\equiv T^\dag/T_\text{slow}^\dag.
\end{eqnarray}
Then the given velocity and the related solution of \eqref{TT-ode} are chosen in the form
\begin{eqnarray}\label{u-x}
\vu^\dag = U^\dag \vu (\vx, s, \tau),\quad \vx^\dag = L^\dag \vx (s, \tau);
\end{eqnarray}
the $\tau$-dependence is $2\pi$-periodic whereas, in general, the $s$-dependence is not periodic.

Transforming Eqn. \eqref{TT-ode} to dimensionless variables and then using the chain rule gives
\begin{eqnarray}
&&{d\vx}/{dt}=\vu,\label{ODE-DLESS}\\
&&\left(\frac{\partial}{\partial\tau}+\mu\frac{\partial}{\partial s}\right)\vx=\varepsilon\vu,\quad \mu\equiv S/\omega =T_\text{fast}/T_\text{slow},\ \varepsilon\equiv{1}/{\omega}.\label{ODE-S}
\end{eqnarray}
The natural small parameter in our consideration is $\varepsilon\equiv 1/\omega$.
The essence of the two-timing method is based on the suggestion that the ratio of two time-scales $\mu$
also represents a small parameter, where $S=S(\omega)$  can be an unknown function of $\omega$.
As a result, Eq.\eqref{ODE-S} contains two independent small parameters, $\mu$ and $\varepsilon$
\begin{eqnarray}\label{eps-PDE}
\vx_\tau+\mu\vx_s=\varepsilon \vu,
\end{eqnarray}
where the subscripts $s$ and $\tau$ denote partial derivatives.
Then, in the two-timing method, we assume that
\begin{eqnarray}
&&\emph{$\tau$ and $s$ are mutually independent variables.}\label{1-key}
\end{eqnarray}
This temporary accepted assumption converts the ODE with two mutually dependent variables $\tau=\tau(t)$ and $s=s(t)$ \eqref{eps-PDE} to an \emph{auxiliary PDE} with two independent variables $\tau$ and $s$ \eqref{eps-PDE}; both look the same but represent  different mathematical objects.
For brevity, we do not discuss any initial data, some their versions can be found in \cite{San, OMalley}.

\section{Distinguished limit procedure applied to vibrodynamics}

Let us consider the PDE \eqref{eps-PDE} with two small parameters $\mu$ and $\varepsilon$.
To construct a rigorous asymptotic procedure for $(\mu, \varepsilon)\to (0,0)$, we have to consider the various paths approaching the origin in the $(\mu, \varepsilon)$-plane (the usual sequence of limits  $\mu\to 0$ and then $\varepsilon\to 0$, or with the order reversed, correspond to the `broken' paths).
A parametrization of a path may be chosen as $\mu= \mu(\varepsilon)$ with the limit $\varepsilon\to 0$ leading to $(\mu, \varepsilon)\to (0,0)$.
In our case, we reduce two dimensionless parameters $\mu, \varepsilon$
to a single one by introducing a parametrization
\begin{eqnarray}\label{eps-PDE-1}
\mu=\varepsilon S(\varepsilon)=\varepsilon^\lambda,\quad\text{where}\ \lambda=\const>0,
\end{eqnarray}
that represents a family of curves in the plane $(\mu, \varepsilon)$; the restriction $\lambda>0$ is required by $\mu\ll 1$.
Eqn.\eqref{eps-PDE-1} means that the slow variable and slow time-scale appear as
\begin{eqnarray}\label{slow-var}
 s=t S(\varepsilon)\equiv t \varepsilon^{\lambda-1},\quad T_\text{slow}\equiv  1/S=\varepsilon^{1-\lambda}
\end{eqnarray}
which shows that the value of $\lambda$ defines the scale of the slow-time variable in terms of the basic small parameter $\varepsilon$.
Then, Eqn.\eqref{eps-PDE-1} transforms \eqref{eps-PDE} to
\begin{eqnarray}
&&{\vx}_{\tau}+ {{\varepsilon^{\lambda}\vx}}_{s} =\varepsilon\vu.\label{LAMBDA-EQN}
\end{eqnarray}
Its solution  $\vx(s,\tau)$ may be formally written as a sum,
\begin{eqnarray}\label{0-series}
\vx(s,\tau)=\sum_{m=0}^M \varepsilon^m \vx_m(s, \tau),\quad m=0,1,2,\dots, M,
\end{eqnarray}
where $M\ge 0$ is an integer; for infinite $M$ \eqref{0-series} represents a regular series, while for a finite $M$ -- a polynomial in $\varepsilon$ or an asymptotic approximation.
%Existence of solution in the form \eqref{0-series} represents a key mathematical assumption.
The Taylor series for  $\vu=\vu(\vx,s,\tau)=(u_{1},u_{2},u_{3})$ is
\begin{eqnarray}\label{0-Taylor}
&&\vu=\vu_0+\varepsilon\,  (\vx_{1}\cdot\nabla_0) \vu_{0}
+\varepsilon^2 \bigl((\vx_{2}\cdot \nabla_{0})\vu_{0} +\frac{1}{2} x_{1k}x_{1l}{\partial^2 \vu_{0}}/{\partial x_{0k}\partial x_{0l}} \bigr)+O(\varepsilon^3),\\
&&\vu_{0}=(u_{01},u_{02},u_{03})\equiv \vu (\vx_0,s,\tau),\quad \vx_0=(x_{01},x_{02},x_{03}), \ \nonumber\\
&&\nabla_0\equiv(\partial/\partial x_{01},\partial/\partial x_{02},\partial/\partial x_{03}),\quad \ f_{0\tau}\equiv\partial f_0/\partial\tau, \ \text{\emph{etc.}},\nonumber
\end{eqnarray}
where we use the summation convention and both vectorial and subscript notations.

In order to make analytic progress, we introduce some common terminology and notation.
First, we accept that any dimensionless (scalar, vectorial, or tensorial) function $g(\vx, s, \tau)$ or $g(s, \tau)$:

$\circ$\ is of order one, $g\sim {O}(1)$; and all its required $\vx$-,  $s$-, and $\tau$-derivatives are also ${O}(1)$;

$\circ$\  is $2\pi$-periodic in $\tau$, \emph{i.e.}  $g(s, \tau)=g(s, \tau+2\pi)$;

$\circ$\  has an average given by
\begin{equation}\label{aver}
\overline{g}\equiv \langle {g}\,\rangle \equiv \frac{1}{2\pi}\int_{\tau_0}^{\tau_0+2\pi}
g(s, \tau)\, d \tau \qquad \forall\ \tau_0;
\end{equation}

$\circ$\  can be split into  averaged  and purely oscillating parts, $g(s, \tau)=\overline{g}(s)+\widetilde{g}(s, \tau)$;
where the \emph{tilde-functions} (or  purely oscillating functions) are such that $\langle \widetilde g\, \rangle =0$ and the \emph{bar-functions} $\overline{g}$ are $\tau$-independent.

Secondly, we introduce the tilde-integration which keeps the result in the tilde-class of functions:
\begin{equation}\label{0-tau-inegr}
\widetilde{g}^\tau\equiv G-\overline{G},\quad G(s,\tau)\equiv\int_0^\tau \widetilde{g}(s,\tau')\, d \tau'.
\end{equation}

Importantly, we impose \emph{five constraints} for the class of solutions we consider:

\noindent
$\bullet$\ All solutions are periodic in $\tau$.

\noindent
$\bullet$\ All solutions can be represented as series or polynomials in $\varepsilon$ \eqref{0-series}.

\noindent
$\bullet$\ The main (zeroth-order) terms $\vu_0$ and $\vx_0$ are not identically zeros or constants:
\begin{eqnarray}\label{0-mainterms}
\vx_0\not\equiv \const,\quad \vu_0\not\equiv \const.
\end{eqnarray}
$\bullet$\ \emph{The amplitude of oscillations (vibrations) of $\vx(s,\tau)$ is small compared with that of averaged motion}:
\begin{eqnarray}\label{0-smal-tilde}
\overline{\vx}_0(s)\not\equiv 0,  \quad \widetilde{\vx}_0\equiv 0.
\end{eqnarray}
$\bullet$\ We introduce a nonvanishing  \emph{drift velocity} ${\overline{\mathbfcal{V}}_2}$
\begin{eqnarray}
&&{\overline{\mathbfcal{V}}_2}={\overline{\mathbfcal{V}}_2}(\overline{\vx}_0(s),s) \equiv \langle(\widetilde{{\vu}}_0^{\tau}\cdot\nabla_0)\widetilde{\vu}_0\rangle=\frac{1}{2}\langle[\widetilde{\vu}_0,\widetilde{\vu}_0^\tau]\rangle\not\equiv 0,\label{drift}
\end{eqnarray}
where $[\vu,\vv]\equiv(\vv \cdot\nabla)\vu-(\vu \cdot\nabla)\vv$ is a commutator of two functions and the subscript `2' at ${\overline{\mathbfcal{V}}_2}$ means that it is quadratic in $\widetilde{\vu}_0$.
\noindent
The requirement \eqref{0-mainterms} is technical; it is intended to keep our  calculations and their results at the lowest order approximations.
The \textbf{\emph{key restriction \eqref{0-smal-tilde} identifies the subject of Vibrodynamics}}, devoted to the studies of large-amplitude average solutions caused by small oscillations/vibrations.
The condition of nonvanishing \emph{drift velocity} \eqref{drift} will be explained in the \emph{Remark 7-1} (\emph{Remark 1} of \emph{Sect.7}).

In the next section we find all available \emph{DL}s by determining all related values of $\lambda$ in \eqref{LAMBDA-EQN} and derive the equations of successive approximations.

\section{Derivations}

Searching for solutions starts by substituting the series \eqref{0-series}
into \eqref{LAMBDA-EQN} whilst taking into account  \eqref{0-Taylor}-\eqref{drift}.
To  obtain the self-consistent successive approximations, one should consider the compatibility condition between the coefficients $\varepsilon$  and $\varepsilon^\lambda$ in \eqref{LAMBDA-EQN}, and $\varepsilon^m$ (with an integer $m$) in \eqref{0-series}.
Namely, a solution in the form \eqref{0-series} is possible only when the exponents $\lambda$ and $1$ in \eqref{LAMBDA-EQN} are multiple of each other.
It means that
\begin{eqnarray}\label{1-lambdaN}
\lambda=N\quad  \text{or}\quad  \lambda=1/N \quad\text{for}\quad N=1,2,3, \dots ,
\end{eqnarray}
where $N=0$ is excluded, since $\lambda>0$.
These two cases simplify \eqref{LAMBDA-EQN} to the equations
\begin{eqnarray}
&\text{Case}\ \lambda=N:\quad & {\vx}_{\tau}+ \varepsilon^N\vx_{s} =\varepsilon\vu,\label{1-PDE1-N}\\
&\text{Case}\ \lambda=1/N:\quad &  {\vx}_{\tau}+ \delta\vx_{s} =\delta^N\vu, \quad \delta\equiv \varepsilon^{1/N},\label{1-PDE1-N-1}
\end{eqnarray}
where for \eqref{1-PDE1-N-1} we must use the series or polynomial \eqref{0-series} in $\delta$, not in $\varepsilon$.
Next, we consider solutions to  \eqref{1-PDE1-N}, \eqref{1-PDE1-N-1} separately for each value $N$:

$\Box$\ \emph{Case} $\lambda=N=1/N=1$. Here we get the equation:
\begin{eqnarray}
&&{\vx}_{\tau}+ \varepsilon\vx_{s} =\varepsilon\vu \label{2-PDE1-A1}
\end{eqnarray}
The substitution of \eqref{0-series}, \eqref{0-Taylor} into \eqref{2-PDE1-A1}
and step-by-step consideration yield
\begin{eqnarray}
&&\overline{\vx}_{0s}=\overline{\vu}_0,\quad \widetilde{\vx}_0\equiv 0\label{2-final-eqn-0}\\
&&\overline{\vx}_{1s}=(\overline{\vx}_1\cdot\nabla_0)\overline{\vu}_0+{\mathbfcal{V}}_2,\quad \widetilde{\vx}_1=\widetilde{\vu}_0^\tau\label{2-final-eqn-1}
\end{eqnarray}
with $\overline{\vu}_0\equiv\langle\vu(\vx_0,s,\tau)\rangle=\langle\vu(\overline{\vx}_0(s),s,\tau)\rangle$.
The derivation of \eqref{2-final-eqn-0},\eqref{2-final-eqn-1} follows the successive steps of orders $\varepsilon^0, \varepsilon^1, \varepsilon^2$:

\emph{Step} $\varepsilon^0$. The zeroth approximation of \eqref{2-PDE1-A1} gives $\widetilde{\vx}_{0\tau}= 0$.
Its $\tau$-integration yields $\widetilde{\vx}_0=0$, where an arbitrary constant vanishes due to zero $\tau$-average.
Hence, we get $\widetilde{\vx}_0\equiv 0$ \eqref{2-final-eqn-0}, and $\vx_0=\overline{\vx}_0(s)$ with an undefined function $\overline{\vx}_0(s)$.

\emph{Step} $\varepsilon^1$. The first approximation yields $\vx_{1\tau}+\vx_{0s}=\vu_0$, where $\vu_0\equiv\vu(\overline{\vx}_0,s,\tau)$.
Its bar-part and tilde-part lead to the first equation in \eqref{2-final-eqn-0} and second equation in \eqref{2-final-eqn-1}.

\emph{Step} $\varepsilon^2$: The second approximation gives $\vx_{2\tau}+\vx_{1s}=(\vx_1\cdot\nabla_0)\vu_0$. Its averaging, together with the results of two previous approximations, lead to the first equation in \eqref{2-final-eqn-1}.

The vibrogenic term ${\overline{\mathbfcal{V}}_2}$ in \eqref{2-final-eqn-1}   represents \emph{the drift velocity} \eqref{drift}.
We have also performed the steps of orders $\varepsilon^3, \varepsilon^4$;  they produce similar but more cumbersome equations that we do not show.
At the same time, the general structure of calculations indicates that the equations of any approximation can be derived.

$\Box$\ \emph{Case} $\lambda=N=2$: Equation \eqref{1-PDE1-N} takes the form:
\begin{eqnarray}
&&{\vx}_{\tau}+ \varepsilon^2\vx_{s} =\varepsilon\vu. \label{2-A2-PDE-0}
\end{eqnarray}

\emph{Step} $\varepsilon^0$. The results here are the same as above: $\widetilde{\vx}_0=0$ and $\vx_0=\overline{\vx}_0(s)$.

\emph{Step} $\varepsilon^1$.  The equation of first approximation is $\vx_{1\tau}=\vu_0$.
Its tilde-part yields $\widetilde{\vx}_1=\widetilde{\vu}_0^\tau$.
The bar-part produces an additional constraint
\begin{eqnarray}\label{bu=00}
&&\overline{\vu}_0\equiv \overline{\vu}_0(\overline{\vx}_0(s),s)\equiv \langle \vu_0(\overline{\vx}_0,s,\tau)\rangle\equiv 0
\end{eqnarray}
which means that the averaged velocity vanishes along the \emph{slow trajectory} $\overline{\vx}_0=\overline{\vx}_0(s)$.
In particular, it means that $\overline{\vu}_{0s}+(\overline{\vx}_{0s}\cdot\nabla_0)\overline{\vu}_0=0$.
One can check that \eqref{bu=00} does NOT allow to obtain a closed system of equations.
We avoid this difficulty by accepting the degeneration of velocity field
\begin{eqnarray}\label{bu=0}
&&\langle\vu(\vx, s,\tau)\rangle_{\vx, s}\equiv 0,
\end{eqnarray}
where the subscripts mean that the average over $\tau$ \eqref{aver} is taken for the fixed variables $\vx,s$;
in particular, it means that
\begin{eqnarray}\label{bu=0-1}
&&\overline{\vu}_{0}\equiv 0,  \  {\partial \overline{\vu}_{0}}/{\partial x_{0k}}\equiv 0, \,  {\partial^2 \overline{\vu}_{0}}/{\partial x_{0k}\partial x_{0l}}\equiv 0, etc.
\end{eqnarray}

\emph{Step} $\varepsilon^2$: The equation of second approximation is $\vx_{2\tau}+\vx_{0s}=(\vx_1\cdot\nabla_0)\vu_0$.
Its bar-part is $\overline{\vx}_{0s}=(\overline{\vx}_1\cdot\nabla_0)\overline{\vu}_0+\langle(\widetilde{\vx}_1\cdot\nabla_0)\widetilde{\vu}_0\rangle$.
The use of results of previous steps yields the full set of equations of zeroth approximation
\begin{eqnarray}
&&\overline{\vx}_{0s}={\overline{\mathbfcal{V}}_2}(\overline{\vx}_{0},s),\quad \widetilde{\vx}_0\equiv 0,\label{B2-final-eqn-1}
\end{eqnarray}
where the fist one represents a closed equation for the averaged motion $\overline{\vx}_0(s)$.
Again, presenting  results of successive approximations has been deliberately stopped after the appearance of the first vibrogenic term ${\overline{\mathbfcal{V}}_2}$, which requires considering only the steps of orders $\varepsilon^0, \varepsilon^1, \varepsilon^2$.
The steps of order $\varepsilon^3$ and above produce similar but more cumbersome equations.
Their algebraic structure indicates that the equations of any approximation can be derived.

$\Box$\ \emph{Case} $\lambda=N$  for $N\ge 3$: Here one can check that the equations of successive approximations are not closed, they do not even allow to obtain a closed equation for $\overline{\vx}_0(s)$.
Therefore all these cases are excluded  from the list of \emph{DL}s.

$\Box$\ \emph{Case} $\lambda=1/N=1/2$  for $N=2$: Equation \eqref{1-PDE1-N-1} takes the form:
\begin{eqnarray}
&&{\vx}_{\tau}+ \delta\vx_{s} =\delta^2\vu.\label{2-A2-PDE-0-}
\end{eqnarray}
As above, the zero-order approximation yields $\widetilde{\vx}_0\equiv 0$.
However, the first-order approximation yields ${\overline{\vx}_{0s}}\equiv 0$ or ${\overline{\vx}_0}\equiv \const$, which contradicts  \eqref{0-mainterms}, \eqref{0-smal-tilde}.
One can also notice that it leads to the degeneration of Taylor series \eqref{0-Taylor} (since all the coefficients are calculated at ${\overline{\vx}_{0}}$) and consequently to the degeneration of all procedures of successive approximations.
Hence $N=2$ is excluded from our study.
Similar outcome takes place for any $N> 2$.
Henceforth, among all cases $\lambda=1/N$ \eqref{1-PDE1-N-1} the only available \emph{DL} is $N=1$, which has been already presented at \eqref{2-final-eqn-0},\eqref{2-final-eqn-1}.

In the next two sections we briefly summarise the results and formulate an instruction for the practical handling of ODEs with different classes of $\vu$.

\section{Summary of results}

\emph{There are only two distinguished limits} \emph{DL-1} and \emph{DL-2} for the auxiliary PDE \eqref{eps-PDE} that is subject to constraints \eqref{0-mainterms}-\eqref{drift}.
We have excluded  all other limits (or all other values of $\lambda$ in \eqref{LAMBDA-EQN}) by regular consideration, so that this result can be seen as a `softly formulated' theorem of uniqueness.
After  \emph{DL-1} and \emph{DL-2} have been identified, both the averaged and oscillatory parts of the equations of successive approximations can be derived.
We have conducted the explicit calculations of averaged equations of the orders $\varepsilon^0, \varepsilon^1, \varepsilon^2$ for \emph{DL-1} and of the orders $\varepsilon^0, \varepsilon^1$ for \emph{DL-2},
while the general structure of the calculations indicates that both the averaged and oscillatory parts of the equations can be derived in any approximation.

\emph{The averaged} \emph{DL-1}\emph{-equations of the first two orders of approximations}
\begin{eqnarray}
&&\overline{\vx}_{0t}=\overline{\vu}(\overline{\vx}_0,t),\label{dl1-aveqn-0}\\
&&\overline{\vx}_{1t}=(\overline{\vx}_1\cdot\nabla_0)\overline{\vu}(\overline{\vx}_0,t)+{\mathbfcal{V}}_2(\overline{\vx}_0,t)\label{dl1-aveqn-1}
\end{eqnarray}
are valid for $\langle\vu(\vx, t,\tau)\rangle_{\vx, t}\not\equiv 0$, see \eqref{0-mainterms},\eqref{0-smal-tilde}.
The unknown functions are $\overline{\vx}_0(t)$ and $\overline{\vx}_1(t)$.
The zeroth approximation \eqref{dl1-aveqn-0} (for the unknown function $\overline{\vx}_0(t)$) formally coincides with the original equation \eqref{ODE-DLESS}.
The averaged equation of the first approximation \eqref{dl1-aveqn-1} is linear (with respect to the unknown function $\overline{\vx}_1(t)$) and contain an `external driving term' ${\overline{\mathbfcal{V}}_2}$ \eqref{drift} that depends only on the previous approximation $\overline{\vx}_0(t)$.
For \emph{DL-1}, Eqs.\eqref{ODE-S},\eqref{eps-PDE-1},\eqref{slow-var} give
\begin{eqnarray}
&&s=t\quad \text{and}\quad T_\text{slow}=1
\label{sT1}
\end{eqnarray}
\emph{The averaged} \emph{DL-2}\emph{-equation of the zeroth-order approximation}
\begin{eqnarray}
&&\overline{\vx}_{0s}={\overline{\mathbfcal{V}}_2}(\overline{\vx}_0,s)\label{dl2-aveqn-0}
\end{eqnarray}
is valid for a degenerated velocity field $\langle\vu(\vx, s,\tau)\rangle_{\vx, s}\equiv 0$ \eqref{bu=0}.
This equation includes ${\overline{\mathbfcal{V}}_2}$ and contains the only unknown function $\overline{\vx}_0(s)$.
Deriving and using such types of equations represents the main aim and the main tool in applications.
For \emph{DL-2}, Eqs.\eqref{ODE-S},\eqref{eps-PDE-1},\eqref{slow-var} give
\begin{eqnarray}
&&s=\varepsilon t=t/\omega\quad \text{and}\quad T_\text{slow}=1/\varepsilon=\omega
\label{sT2}
\end{eqnarray}

The calculations of next successive approximations can be successfully continued.
The unknown function in the averaged equation of $n$-th approximation ($n=0,1,2,3,\dots$) is $\overline{\vx}_n(s)$, while all equations for $n=1,2,3,\dots$ are of the first-order and linear.
These equations consist of \emph{two universal groups of terms}:

\emph{Group 1:}  All the terms formally coinciding with the small-amplitude decompositions of the original equation \eqref{ODE-DLESS} obtained in the absence of any oscillations;
these terms contain  the averaged unknown functions  $\overline{\vx}_0(s),\overline{\vx}_1(s),\overline{\vx}_2(s)$, \emph{etc.};

\emph{Group 2:} Various universally defined vibrogenic terms, where the first one by order of appearance is ${\overline{\mathbfcal{V}}_2}={\overline{\mathbfcal{V}}_2}(\overline{\vx}_0(s),s)$.

\noindent
The simplest examples of such \emph{universal structures of averaged equations} are given by  Eqs. \eqref{dl1-aveqn-0}-\eqref{sT2}.

There are three \emph{key remarks to clarify the results:}

\emph{Remark 5-1:} The mathematical justifications of the used method are given in the papers and books, quoted in the \emph{Introduction.}
In this paper, we \emph{formally} (without error estimates and related mathematical formulations) accept the validity of all derived approximations and equations.

\emph{Remark 5-2:} In particular, we accept the validity of the inhomogeneous linear  Eqn. \eqref{dl1-aveqn-1} with ${\overline{\mathbfcal{V}}_2}$ as an `external driving term'.
Alternatively, one may perform a standard amplitude linearization of \eqref{dl1-aveqn-0} that gives only the homogeneous version of Eqn. \eqref{dl1-aveqn-1} with ${\overline{\mathbfcal{V}}_2}\equiv 0$.
Therefore such a standard linearization is unnecessary: all related perturbations are already described by  \eqref{dl1-aveqn-1} as its homogeneous solutions.
At the same time, any particular solution of the inhomogeneous Eqn.\eqref{dl1-aveqn-1} give us an additional linear perturbation that cannot be obtained by any linearization of Eqn. \eqref{dl1-aveqn-0}.
One can see that Eqn.\eqref{dl1-aveqn-1} leads to a new (and to more general and rigorous) procedure of linearization and to new results for linear perturbations.

\emph{Remark 5-3:}
Comparison between the orders of physical amplitudes of solutions to \eqref{dl1-aveqn-0}, \eqref{dl1-aveqn-1} and \eqref{dl2-aveqn-0}
can be made by returning Eqn.\eqref{dl2-aveqn-0} in the original form of Eqn.\eqref{ODE-DLESS} rewritten with the use of physical variable $t$.
Indeed, any equation $d\,\overline{x}/ds=\overline{u}$ with $s\equiv\varepsilon t$ can be rewritten as $d\,\overline{x}/dt=\varepsilon \overline{u}$.
This shows that an equation $d\,\overline{x}/ds=\overline{u}$ of the order $O(1)$ actually describes a solution $\overline{x}(t)$ of order $O(\varepsilon)$.
Hence, Eqn.\eqref{dl2-aveqn-0} and Eqn.\eqref{dl1-aveqn-1}  describe the amplitudes of the same order $O(\varepsilon)$, despite Eqn.\eqref{dl2-aveqn-0} formally exhibiting order $O(1)$.

\section{Practical instruction}

The practical outcome is the ability to deal with a given dimensionless ODE $d\vx/dt=\vu$ \eqref{ODE-DLESS} containing a particular oscillating velocity $\vu=\vu(\vx,s,\tau)$,  where $s=t/\omega^\nu (\nu=\const>-1)$ and $\tau\equiv \omega t$. The equations can be solved as $\vx=\vx(s,\tau)$ only for $\nu=0$ and $\nu=1$, for all others values of $\nu$ the solutions (that are subject of \eqref{0-mainterms}-\eqref{drift}) do not exist.
\emph{The explicit instruction is}:

$\bullet$\textbf{ \emph{Case A}:}
If an equation belongs to the class
\begin{eqnarray}\label{results1}
&& d\vx/dt=\vu(\vx,t,\tau), \quad\text{with}\quad\langle{\vu}\rangle_{\vx,t}\not\equiv 0,
\end{eqnarray}
then one has to use \emph{\text{DL-1}}.
The subscripts at $\langle\dots\rangle$ mean that the average over $\tau$ \eqref{aver} is taken for the fixed variables $\vx,t$.
The related closed system of averaged equations is \eqref{dl1-aveqn-0}-\eqref{sT1}.

$\bullet$ \textbf{\emph{Case B}:}
If  an equation belongs to the class with a degeneration of $\vu$ as
\begin{eqnarray}
&& d\vx/dt=\vu(\vx,t/\omega,\tau),\quad  \text{with}\quad \langle{\vu}\rangle_{\vx,t}\equiv 0,\label{results2}
\end{eqnarray}
then one has to use \emph{\text{DL-2}}.
The closed system of averaged equations \eqref{dl2-aveqn-0},\eqref{sT2} appears already in the zeroth approximation.
The argument $t/\omega$ may look artificial.
A more useful case could be a purely oscillating velocity, as in the following

$\bullet$ \textbf{\emph{Case C}:}
If  an equation belongs to the class with another degeneration of $\vu$ as
\begin{eqnarray}
&& d\vx/dt=\vu(\vx,\tau),\label{results3}
\end{eqnarray}
then one has to use \emph{DL-1} for $\langle{\vu}\rangle_{\vx}\not\equiv 0$ or \emph{DL-2} for $\langle{\vu}\rangle_{\vx}\equiv 0$, with the same averaged equations as in \emph{Cases A and B}.

In practice, one has to produce the dimensionless form of a velocity field $\vu(\vx,s,\tau)$ that appears in a chosen particular problem,
and then to solve the appropriately chosen \emph{DL-1} or \emph{DL-2} averaged equations which lead to particular expressions for $\overline{\vx}_0(t)$ and $\overline{\vx}_1(t)$.
After that, one can take into account the boundary conditions and switch to the dimensional ODE \eqref{TT-ode}; those steps are not considered here.

\section{Drift velocity}

\emph{The  drift velocity} $\overline{\mathbfcal{V}}_2=\overline{\mathbfcal{V}}_2(\overline{\vx}_0(s),s)$ \eqref{drift} plays a key part in both theory and  experiments.
In our study it appears as the first vibrogenic (averaged and nonlinear in $\widetilde{\vu}_0$) term of the averaged equations.
One can check that the accepted condition $\overline{\mathbfcal{V}}_2\not\equiv 0$ \eqref{drift} is valid for the vast majority of prescribed fields $\widetilde{\vu}_0$
(an exception is the separating of $\tau$-variable, such as $\vu=f(\tau)\vv(\vx,s)$, where $\overline{\mathbfcal{V}}_2\equiv 0$; related physical examples are standing waves).

\emph{The universality of} ${\overline{\mathbfcal{V}}_2}$ exhibits itself as:

(i) it is defined in \eqref{drift} for a general analytic expression of oscillating velocity $\widetilde{\vu}_0$, hence many particular derivations of averaged equations in vibrodynamics, see \emph{e.g.} \newline \cite{Blekhman1, Thom, Blekhman2}, may be replaced or optimized by the use of this formula;

(ii) this expression is the same for both \emph{DL-1} and \emph{DL-2}; for velocity fields with $\overline{\mathbfcal{V}}_2\equiv 0$ it leads to new \emph{DLs} and to the drift velocities of higher orders (that is also expressed by universal formulae), see \emph{Remark 7-1} below;

(iii) $\overline{\mathbfcal{V}}_2$ appears as an universal term in the universally defined averaged equations;

(iv) it is applicable for any fixed spatial domain $\mathcal{D}$ of definition of $\vu(\vx,t)$;
it is remarkable that a no-leak condition $\vu\cdot\vn=0$ at the boundary $\partial\mathcal{D}$ (with a normal vector $\vn$) leads to ${\overline{\mathbfcal{V}}_2}\cdot\vn=0$ that does not allow the averaged trajectories to cross the boundary; the generalization of this property for an oscillating boundary is given in \cite{VladProc};

(v) The expression for \eqref{drift} $\overline{\mathbfcal{V}}_2$ contains, for example, the classical Stokes drift generated by periodic water waves, that represent a classical subject in fluid dynamics, magneto-hydrodynamics and plasma physics, see \emph{e.g.} \cite{Stokes, Lamb, Batchelor, Craik, Eames, VladimirovD, Vladimirov1, Brem, Moni}.
%\newline \cite{}.

The key remarks, clarifying \emph{the multiple roles of the drift velocity} are:

\emph{Remark 7-1:} The presence of two distinguished limits \emph{DL-1} and \emph{DL-2} is due to the \emph{enforced condition} ${\overline{\mathbfcal{V}}_2}\not\equiv 0$ \eqref{drift}.
If we allow further degeneration of velocity (when both $\overline{\vu}_0\equiv 0$ and $\overline{\mathbfcal{V}}_2\equiv 0$), then there is a new distinguished limit \emph{DL-3}, corresponding to $\lambda=3$ in \eqref{LAMBDA-EQN}, where \eqref{dl2-aveqn-0}  is replaced by
\begin{eqnarray}
&&\overline{\vx}_{0s}={\overline{\mathbfcal{V}}_3},\quad s=\varepsilon^2 t, \quad \overline{\mathbfcal{V}}_3\equiv\frac{1}{3}\langle[[\widetilde{\vu}_0,\widetilde{\vu}_0^\tau],\widetilde{\vu}_0^\tau]\rangle,\label{V3}
\end{eqnarray}
where a new drift velocity $\overline{\mathbfcal{V}}_3$ is cubic in $\widetilde{\vu}_0$.
Such successive degenerations of $\vu$ % ($\overline{\vu}_0\equiv 0$, ${\overline{\mathbfcal{V}}_2}\equiv 0$, ${\overline{\mathbfcal{V}}_3}\equiv 0$, \emph{etc.})
have been studied by \cite{VladimirovD}, \cite{Vladimirov1}, where a  hyperbolic PDE was considered
\begin{eqnarray}\label{a-eqn}
\partial a/\partial {t^\dag}+(\vu^\dag\cdot\nabla^\dag) a=0
\end{eqnarray}
describing the transport of a passive scalar admixture (or a lagrangian marker) $a(\vx^\dag, t^\dag)$ by the same velocity field as in \eqref{TT-ode},\eqref{scales}.
The ODE \eqref{TT-ode} describes the characteristic curves for \eqref{a-eqn}, hence the appearance of the same drift velocities in both problems can be expected.
Such a procedure of successive degenerations of velocity produces a sequence of distinguished limits \emph{DL-N} with $N=1,2,3,\dots$ and with different drift velocities.

\emph{Remark 7-2:} We have established that \emph{DL-1} is usable for $\overline{\vu}_0\not\equiv 0$, \emph{DL-2} -- for $\overline{\vu}_0\equiv 0$ but $\overline{\mathbfcal{V}}_2\not \equiv 0$ and \emph{DL-3} -- for  $\overline{\vu}_0\equiv 0$, $\overline{\mathbfcal{V}}_2\equiv 0$ but $\overline{\mathbfcal{V}}_3\not\equiv 0$.
The \emph{physical meaning}  can be clarified  by comparison between the displacement $\Delta\overline{\vx}$ of a particle during the time-interval $\Delta t=1$:

(i) for a flow with $\overline{\vu}_0\not\equiv 0$ the slow time $s=t$ shows that a particle driven by $\overline{\vu}_0\sim O(1)$ is displaced by $\Delta\overline{\vx}\sim O(1)$;.

(ii) for a flow with $\overline{\vu}_0\equiv 0$ and ${\overline{\mathbfcal{V}}_2}\not\equiv 0$ the slow time $s=\varepsilon t$ means that a particle driven by the drift velocity ${\overline{\mathbfcal{V}}_2}\sim O(1)$,
is displaced by $\Delta\overline{\vx}\sim O(\varepsilon)$;

(iii) for a flow with  $\overline{\vu}_0\equiv 0$,  ${\overline{\mathbfcal{V}}_2}\equiv 0$ and ${\overline{\mathbfcal{V}}_3}\not\equiv 0$ the slow time $s=\varepsilon^2 t$ means that a particle driven by the drift velocity ${\overline{\mathbfcal{V}}_3}\sim O(1)$ \eqref{V3} is displaced by $\Delta\overline{\vx}\sim O(\varepsilon^2)$.

\noindent
Such interrelations between the amplitudes of solutions $\overline{\vx}(s(t))$ described by  zero-order equations with different $s(t)$ have also been clarified in \emph{Remark 5-3}.

\emph{Remark 7-3:} There is a linear in $s$ growth of $\overline{\vx}_1(s)$  within \emph{DL-1} by virtue of \eqref{dl1-aveqn-1} due to its `external driving term' $\overline{\mathbfcal{V}}_2$, at least for some particular functions $\widetilde{\vu}_0$.
Such a growth can be called as the \emph{drift-initiated} one.
Remarkably,  such a growth does not have a character of instability, it appears as an externally driven motion.
For example, if the solution of the main approximation possesses a fixed point  with $\overline{\vx}_{0}(s)=\overline{\vx}_{0}^*\equiv \const$, then, for  ${\overline{\mathbfcal{V}}_2}(\overline{\vx}_{0}^*,s)\neq 0$, a linear in $s$ growth $\overline{\vx}_1\sim s\,{\overline{\mathbfcal{V}}_2}$ follows immediately.
In general, the presence of $\overline{\mathbfcal{V}}_2 \not\equiv 0$ deforms the structure of zeroth-order trajectories.
One may expect that such a linear growth can have a practical meaning.

\emph{Remark 7-4:} The setting of \emph{ general notions of stability or instability of averaged motion and averaged equilibria} was considered by \cite{Bog} in the framework of their method.
However, as it has been shown for \emph{DL-1}, even the existence of an averaged equilibria in zero approximation $\overline{\vx}_0(s)=\overline{\vx}_{0}^*\equiv \const$ is compromised by the presence of drift velocity, leading to a linear growth of $\overline{\vx}_1(s)$ \eqref{dl1-aveqn-1}.
A similar result \emph{for DL-2} is also valid, but not shown here.
Hence, the notion of stability is not adapted to our consideration yet.

\subsection{Examples of drifts}

We have adopted the  terminology of a moving particle in a  three-dimensional velocity field for the sake of brevity and convenience only.
In general, the system of equations \eqref{ODE-DLESS}  can have any applied meaning and any dimension.
Certainly, the term \emph{`drift'} itself can be confusing when it is not linked to any physical motion.
Let us consider the drifts in two most popular model of mathematical biology (see \cite{Murray}).
In biological applications, the time-oscillation of related coefficients can be caused, for example, by the day-night or seasonal variations of temperature.

\noindent
\emph{Example 1: The Logistic equation} is

\begin{eqnarray}\label{logistic}
\frac{dx}{dt}=ax(1-bx)=ax-ab x^2,
\end{eqnarray}
where $a,b$ are experimentally defined constants.
Eqn.\eqref{logistic} represents a one-dimensional case of \eqref{ODE-DLESS} with the only component of `velocity':
\begin{eqnarray}\label{logistic-u}
u(x,t)=ax-ab x^2.\quad
\end{eqnarray}
In our study we take both $a\sim O(1)$ and $b\sim O(1)$ as $2\pi$-periodic functions of $\tau$
\begin{eqnarray}\label{ab}
 a=\overline{a}+\widetilde{a}(\tau),\ b=\overline{b}+\widetilde{b}(\tau);\quad \overline{a}=\const,\ \overline{b}=\const,
\end{eqnarray}
and $x(s, \tau)$ as a one-dimensional series \eqref{0-series}.
The asymptotic equations are given by \emph{DL-1} and \emph{DL-2}.
For the latter, it must be $\overline{a}=\overline{b}=0$ and $\widetilde{a}\sim \widetilde{b}\sim O(1)$.
The calculation of the drift velocity \eqref{drift} gives
\begin{eqnarray}\label{drift-logis}
{\overline{\mathcal{V}}_2} =-K \overline{x}_0^2;\quad K\equiv \langle \widetilde{a}^\tau \widetilde{a}\, \widetilde{b}\,\rangle-\overline{a}\langle \widetilde{a} \widetilde{b}^\tau\rangle,
\end{eqnarray}
where $\overline{x}_0(s)$ and $\overline{x}_1(s)$ satisfy \eqref{dl1-aveqn-0},\eqref{dl1-aveqn-1} or \eqref{dl2-aveqn-0}
with different $s$ for \emph{DL-1} or \emph{DL-2}, in accordance with \eqref{results1}-\eqref{results3}.
The \emph{DL-2} produces an averaged ODE
\begin{eqnarray}\label{PPeqn2}
d\,\overline{x}_0/ds=-K\, \overline{x}_0^2\quad\text{where}\quad K\equiv \langle \widetilde{a}^\tau \widetilde{a}\, \widetilde{b}\,\rangle,\quad s=\varepsilon t,\quad \overline{x}_0(0)>0
\end{eqnarray}
Its solution
\begin{eqnarray}\label{PPeqnSol2}
\overline{x}_0(s)=\overline{x}_0(0)/(1+\overline{x}_0(0) K s)
\end{eqnarray}
is monotonically decreasing for $K>0$ and increasing for $K<0$.
The latter shows the increasing of population $\overline{x}_0(s)$ in the biologically challenging case of purely oscillating coefficients.
The \emph{DL-1} produces averaged  ODEs
\begin{eqnarray}\label{PPeqn1}
&&d\,\overline{x}_0/dt=\overline{a}\, \overline{x}_0-\overline{c}\,\overline{x}_0^2 \quad\text{where}\quad \overline{a}>0,\quad s= t,\quad \overline{x}_0(0)>0\\
&&d\, \overline{x}_1/dt=\overline{x}_1(\overline{a}-2\overline{c}\,\overline{x}_0)-K\, \overline{x}_0^2\quad\text{where}\quad\overline{x}_1(0)=0\nonumber
\end{eqnarray}
where $c\equiv ab$, $\overline{c}=\overline{a}\,\overline{b}+\langle \widetilde{a}\,\widetilde{b}\rangle$.
The first equation gives the well-known analytic solution $\overline{x}_0(t)$ that rapidly approaches a constant for $t>\overline{a}$, see \cite{Murray}.
Then the second equation produces an exponential homogeneous solution (related to the standard properties of stability or instability) and the linear in $t$ drift-initiated growth as a particular solution, mentioned in \emph{Remark 7-3}.

\vskip 1mm
\noindent
\emph{Example 2: The Predator-Prey equation} is
\begin{eqnarray}\label{PPeqn}
  && dx/dt=\alpha x-\beta xy \\
  && dy/dt=-\gamma y +\mu xy,
\end{eqnarray}
where $\alpha=\overline{\alpha}+\widetilde{\alpha}$, $\beta=\overline{\beta}+\widetilde{\beta}$, $\gamma=\overline{\gamma}+\widetilde{\gamma}$ and $\mu=\overline{\mu}+\widetilde{\mu}$
with constants $\overline{\alpha}, \overline{\beta}, \overline{\gamma}, \overline{\mu}$ all of order one for \emph{DL-1} or all zero for \emph{DL-2}.
The main term of the oscillating part of `velocity' is:

\begin{eqnarray}\label{u-pred}
\widetilde{\vu}_0=\left(
                   \begin{array}{c}
                      \widetilde{\alpha} \overline{x}_0-\widetilde{\beta} \overline{x}_0 \overline{y}_0\\
                     -\widetilde{\gamma} \overline{y}_0+\widetilde{\mu} \overline{x}_0 \overline{y}_0\\
                   \end{array}
                 \right).
\end{eqnarray}
It gives the drift velocity
\begin{eqnarray}\label{drift-pred}
{\overline{\mathbfcal{V}}_2}=\left(
                   \begin{array}{c}
                     A\overline{x}_0 \overline{y}_0 -B\overline{x}_0^2 \overline{y}_0 \\
                     -C\overline{x}_0 \overline{y}_0+B\overline{x}_0 \overline{y}_0^2\\
                   \end{array}
                 \right),
\end{eqnarray}
where
\begin{eqnarray}
A\equiv\langle\widetilde{\beta}\widetilde{\gamma}^\tau\rangle,\quad B\equiv\langle\widetilde{\beta}\widetilde{\mu}^\tau\rangle,\quad C\equiv\langle\widetilde{\alpha}\widetilde{\mu}^\tau\rangle.
\end{eqnarray}
\noindent
It can be seen that the biological restrictions, that the unknown functions (representing populations) should be non-negative for any positive initial data, are automatically satisfied for both \emph{DL-1} and \emph{DL-2} equations.
%The results like  $\overline{x}_0=0$ or $\overline{y}_0=0$, at a certain instant of time, may lead to the conclusion on total population death or on the lack of generality of presented models.
The equations \eqref{dl1-aveqn-0} and \eqref{dl2-aveqn-0} can be solved analytically, all solutions with positive initial data belong to the first octant $(\overline{x}_0>0, \overline{y}_0>0)$.
For \emph{DL-2}, Eqn.\eqref{dl2-aveqn-0} has an integral $(\overline{x}_0-A/B)\,(\overline{y}_0-C/B)=\const$, ($B\neq 0$); the related
trajectories in the $(\overline{x}_0,\overline{y}_0)$-plane represent hyperbolas  with an unstable equilibrium at $(A/B,C/B)$.
%One can see that the averaged equations \eqref{dl1-aveqn-0}-\eqref{dl2-aveqn-0} in both examples are mathematically tractable (and belonging to the degenerated \emph{Case C} \eqref{results3}), while the level of their nonlinearity can increase.
Some solutions are growing with time, others are decaying.
For \emph{DL-1}, an interesting open question is: should the linearly growing solutions to Eqn.\eqref{dl1-aveqn-1} be considered as biologically meaningful?
However, \emph{the aim of both examples is only to expand the terminological and topical scope of the paper, without going into biological issues.}
More elaborated example from mathematical biology is given by \cite{Mor}, who derived the averaged equations (which include homogenisation) and analyzed their solutions.
In general,  the number of applications may be significant.
For example, several applied systems can be taken from \cite{Stro, Murray}, after installing the time-oscillations into the coefficients.
%The only place that requires   special attention is the scaling \eqref{scales} that can vary for different applications.
Multiple examples of two-timing equations and drifts in fluid dynamics are given by \cite{VladimirovL, VladimirovR1, VladimirovR2, VladProc}.
In particular,  \cite{VladimirovR1}, \newline \cite{VladimirovR2}  show that the self-propulsion of deformed bodies in micro- hydrodynamics represents a drift motion in a self-generated oscillating field.

\section{Discussion}

\emph{Remark 8-1:} There is a \emph{curious mismatch} between the methods of solving two classes of two-timing problems described at the beginning of \emph{Introduction}.
The major tool for finding approximations for nonlinear/modulated oscillations is the forcing of secular terms to vanish, see \emph{e.g.} \cite{San, OMalley}.
In contrary, the secular terms appear in vibrodynamics due to the drift $\overline{\mathbfcal{V}}_2$, which is essentially nonvanishing in the vast majority of flows and plays a key role in fluid mechanics and plasma physics.

\emph{Remark 8-2:} The constraint \eqref{0-mainterms},\eqref{0-smal-tilde} of a \emph{non-vanishing zeroth-order term} can be abolished; if one takes $\overline{\vx}_0\equiv 0$ then similar to above results appear in higher approximations.

\emph{Remark 8-3:} One can introduce in \eqref{ODE-DLESS} a first correction to the given velocity as $\vu(\vx,t)+\varepsilon\vv(\vx,t)$, where $\vv$ is another given function with similar to $\vu$ properties.
Then all previous results  remain principally the same, except for the appearance of additional terms in the averaged equations.
Such additions can be very useful in applications, see \cite{Vladimirov, Yudovich, VladNewt}.

\emph{Remark 8-4:} The theory of this paper does not include the most important class of second-order ODEs with given oscillating forces.
Indeed, for the Newton equation, \eqref{eps-PDE} is replaced with
\begin{eqnarray}
\left(\frac{\partial}{\partial\tau}+\mu\frac{\partial}{\partial s}\right)^2\vx=\varepsilon\vf,\label{second-1}
\end{eqnarray}
where $\vf=\vf(\vx,s,\tau)$ is a given oscillating dimensionless force, replacing the dimensionless velocity $\vu$; small parameters $\mu,\varepsilon$ have the same meaning as in \eqref{eps-PDE}, see \cite{Vladimirov, Yudovich, VladNewt}.
To  represent Eqn. \eqref{second-1} as a system of first-order equations of doubled dimension, we introduce an auxiliary unknown function $\vy=\vy(s, \tau)$.
Then
\begin{eqnarray}
&&\left(\frac{\partial}{\partial\tau}+\mu\frac{\partial}{\partial s}\right)\vx=\vy\label{second-2}\\
&&\left(\frac{\partial}{\partial\tau}+\mu\frac{\partial}{\partial s}\right)\vy=\varepsilon\vf(\vx,s,\tau).\nonumber
\end{eqnarray}
One can see that the six-dimensional system of equations \eqref{second-2} does not belong to the same class as \eqref{eps-PDE}.
Namely, the six-dimensional `generalized velocity' $(\vy, \varepsilon\vf)$ in \eqref{second-2} is of order $O(1)$ for the first three equations and of order $O(\varepsilon)$ for other equations,
while in \eqref{eps-PDE} all components of velocity are of order $O(\varepsilon)$.
In addition, the functional forms of all six  components of `generalized velocity' are strongly degenerated.
This is why the drift velocity ${\overline{\mathbfcal{V}}_2}$ does not appear in the averaged second-order equations with an oscillating force;
it is replaced with a \emph{vibrogenic force} ${\overline{\mathbfcal{F}}_2}$, see \cite{Vladimirov, Yudovich},
\cite{VladNewt}.

\begin{acknowledgments}
\emph{Acknowledgements:}
The author  is grateful to Profs. A.D.D.Craik, I. Eltayeb, K.I.Ilin, D.W.Hughes, D. Kapanadze, H.K.Moffatt, M.T.Montgomery, A.B.Morgulis, T.J. Pedley, M.R.E.Proctor, A.I. Shnilerman and V.I.Yudovich for discussions, and to Mr. A.A.Aldrick for help with the manuscript.
\end{acknowledgments}

\end{document}